\documentclass[12pt]{amsart}

\usepackage{amsmath,amsthm,amsfonts,amscd,amssymb}
\usepackage{times,euler,eucal,eufrak}

\usepackage{hyperref}

\title[Jacobi ensembles and free probability]
{Product of random projections, Jacobi ensembles
and universality problems arising from free probability}
\author {Beno\^\i{}t Collins}
\thanks{B.C. is currently a JSPS postdoctoral fellow}
\address{
Department of Mathematics, Graduate school of Science
Kyoto University, Oiwake-cho, Kitashirakawa, Sakyo-ku
Kyoto 606-8502 - Japan} \email{collins@math.kyoto-u.ac.jp}

\theoremstyle{plain}
\newtheorem{lemma}{Lemma}[section]
\newtheorem{theorem}[lemma]{Theorem}
\newtheorem*{theorem*}{Theorem}
\newtheorem{proposition}[lemma]{Proposition}
\newtheorem{corollary}[lemma]{Corollary}

\theoremstyle{definition}

\newtheorem{assumption}{Assumption}

\theoremstyle{remark}
\newtheorem*{remark}{Remark}

 \DeclareMathOperator{\Tr}{Tr}

\newcommand{\M}[1]{M_{#1}(\mathbb{C})}

\begin{document}

\begin{abstract}
We consider the product of two independent randomly rotated projectors.
The square of its radial part turns out to be distributed as
a Jacobi ensemble. 
We study its global and local properties in the large dimension
scaling relevant to free probability theory. We establish
asymptotics for one point and two point correlation functions,
as well as properties of largest and smallest eigenvalues.
\end{abstract}

\maketitle

\section{Introduction.}

In this paper, we consider the asymptotic distribution of eigenvalues of a
random matrix of the form $\pi_n\tilde\pi_n\pi_n$ where $\pi_n$ and
$\tilde\pi_n$ are independent $n\times n$
random orthogonal  projections, of ranks $q_n$ and $\tilde q_n$, whose
distributions are invariant under unitary conjugation. This question is part of
a more general problem in free probability theory,
where one would like to study matrices of the form
$\pi_n A_n\pi_n$ where $A_n$ is a random matrix whose distribution 
is unitarily invariant, and whose empirical eigenvalues distribution converges.
Indeed, the contraction of a subalgebra by a free projection
has been much studied, and
the pair $(\pi_n, A_n)$ is the most natural asymptotic model 
of a random variable $A_n$ free from a projector $\pi_n$.

Our approach relies on the fact that we can explicitly compute the eigenvalue
distribution of the above model. Similar computations have been initiated in the
paper of Olshanski \cite{0724.22020}
(see also the author's PhD thesis \cite{Co}) but the method presented in this paper is 
more elementary.

Actually we will see that the random matrix
$\pi_n\tilde\pi_n\pi_n$
is distributed according to a Jacobi ensemble of parameters 
$(q_n,n-\tilde q_n-q_n,\tilde q_n -q_n)$.
For the definition, see Equation \eqref{defjac} and for a good review, \cite{Fo}.
We use asymptotic properties of Jacobi polynomials in order to derive the
asymptotic distributions of eigenvalues. We find that the one point function
has an explicit limit, which we relate to free probability theory. We
also check that the universality conjectures of Mehta 
are verified for this
model both in the bulk of the spectrum and at the soft and hard edges
(see \cite{MR92f:82002}, 
conjectures 1.2.1 and 1.2.2 for a statement of this universality problem, and
recent works of \cite{MR2002j:15024,MR2001i:82037} for important 
breakthroughs towards these conjectures).

The universality conjectures at the hard edge in different frameworks have been
established by Kuijlaars and Vanlessen \cite{math-ph/0204006} and our 
result extends a part of their work without using Riemann-Hilbert methods.

As for universality conjectures at the soft edge, 
a recent work of Ledoux \cite{Le} 
gives explicit non asymptotic bound for the tail of the distribution of the 
largest eigenvalue of a modified Jacobi ensemble. 
Our cornerstone result is Theorem \ref{jacobi}:

\begin{theorem*}
Let $X,X'\in\M{q_n}$ be independent Wishart matrices 
of parameters $(q_n,n-\tilde q_n,1/q_n)$ and $(q_n,\tilde q_n,1/q_n)$. 
Let $J=(X+X')^{-1/2}X(X+X')^{-1/2}$ (this is well defined by Lemma \ref{lem2}).
On the other hand, 
let $\pi_n\in\M{n}$ be a constant orthogonal projection of rank $q_n$ and
$\tilde\pi_n\in\M{n}$ be a random uniformly distributed orthogonal projection of
rank $\tilde q_n$.

Then, under the isomorphism $\pi_n\M{n}\pi_n=\M{q_n}$, 
the following equality holds in distribution:
$$\pi_n\tilde\pi_n\pi_n\overset{\mathcal{L}}{=}J$$
In particular for $\tilde q_n\geq q_n$ and $q_n+\tilde q_n \leq n$,
the distribution of  $\pi_n\tilde{\pi}_n\pi_n$ is a
 Jacobi ensemble of parameter $(q_n,n-q_n-\tilde q_n, \tilde q_n -q_n)$ 
 (on $\M{q_n}$).
\end{theorem*}

There is a striking analogy between this result and that of \cite{CC} 
stating results of asymptotic freeness for so-called ``Beta Matrices" 
whose eigenvalue distribution actually follows Jacobi ensembles.
Theorem \ref{jacobi} can also be found under a different formulation and for
different purposes in \cite{Do}.

To the knowledge of the author, the link between products of randomly rotated
projections and Jacobi ensembles had only been observed asymptotically
so far, and not at the finite dimension level.

In accordance to Theorem \ref{jacobi}, we consider Jacobi ensembles of
type $J(n,a_n,b_n)$ with $n\rightarrow\infty$, and let $\Lambda_n^{a_n,b_n}$ 
be the random set of its eigenvalues. 
This random set is a so-called {\it determinantal point process \it}. We call
$K_n^{a_n.b_n}$ the kernel that drives it (see section \ref{defker}).
If
$a_n\sim \alpha n$, $b_n\sim \beta n$, $n\rightarrow\infty$, 
free probabilistic arguments show
that the associated mean counting probability measure converges 
in moments as $n\rightarrow \infty$.

Our first series of results are 
Theorem \ref{onepoint}, Proposition \ref{norme}, Theorem \ref{univ-local},
which we summarize here:

\begin{theorem*}
Assume $a_n\sim \alpha n$, $b_n\sim \beta n$ as $n\to\infty$ 
(assumption \ref{hypo})
\begin{itemize}
\item
The density of the expectation of the eigenvalues counting measure 
for eigenvalues of $J(n,a_n,b_n)$
converges towards the $n=\infty$ limit of the densities 
defined in Equation \eqref{deff}.
This convergence is uniform on any compact set not containing the boundary 
points $r,s$  of the spectrum
(see Equation \eqref{defrs} for the definition of $r,s$).
\item
Let $K_n^{a_n.b_n}$
be the kernel associated to the Jacobi ensemble as a determinantal point process
(for definitions, see section \ref{detpp})
Then, as $n\rightarrow\infty$ and uniformly for $x\in [r+\varepsilon,s-\varepsilon ]$,
($\varepsilon >0$)
and $u,v$ on compact sets,
\begin{equation*}
\frac{1}{nf(x)}K_n^{a_n,b_n}(x+\frac{u}{nf(x)},x+\frac{v}{nf(x)})=
\frac{\sin \pi (u-v)}{\pi(u-v)}+O(n^{-1})
\end{equation*}
where $f$ is defined at Equation \eqref{deff}.
In  other words, the universality conjecture of Mehta holds in the bulk of the
spectrum.
\item 
For any compact set $K$ such that $K\cap [r,s]= \emptyset$, there exists a constant
$C>0$ such that for all $n$,
$P (\Lambda_n^{a_n,b_n}\cap K\neq \emptyset)<e^{-Cn}$ 
\end{itemize}
\end{theorem*}

The limit distribution defined in Equation \eqref{deff} admits a connected
spectrum $[r,s]$ (plus possibly up to two atoms). In addition, the non atomic 
part admits a continuous density that behaves either like $(x-r)^{1/2}$ or $(x-r)^{-1/2}$
close to the spectrum. 
Following conventions in the physics litterature, the first case shall be
 referred as a ``soft edge'' and the latter one
as a ``hard edge''.
We obtain that the relevant spacings for obtaining kernels are the usual
ones ($n^{-2/3}$ for the soft edge, and $n^{-2}$ for hard edge).
Our main theorems are Theorems \ref{airy} and \ref{bessel}:

\begin{theorem*}
\begin{itemize}
\item
At the soft edge, under Assumption \ref{soft}, let $s_n$ be as Equation \eqref{defrsn}
and 
$$h_{n}=\left(
\frac {\sqrt { \left( 1+\alpha_n \right)  \left( 1+\beta_n \right) 
\left( 1+\alpha_n+\beta_n \right) }}{2 \left( 1- s_n^2 \right) ^{2}}\right)^{1/3}$$
Then for any $\varepsilon >0$, one has
\begin{equation*}
\frac{1}{h_{n}n^{2/3}}K_{n}^{a_n,b_n}\left(x+\frac{u}{h_{n}n^{2/3}},x
+\frac{v}{h_{n}n^{2/3}}\right)
=Ai(u,v)+0(n^{-1/3+\varepsilon})
\end{equation*}
where $Ai$ is defined at Equation \eqref{defairy}.
\item
At the hard edge, under Assumption \ref{hard}
(without loss of generality we assume that $r=-1$), 
for any $u,v\in\mathbb{R}_+$,
\begin{equation*}
\frac{1}{2n^2(1+\alpha_n)}K_n^{a_n,b_n}\left(-1+\frac{u}{2n^2(1+\alpha_n)},-1+\frac{v}{2n^2(1+\alpha_n)}\right)
= F_b(u,v)+O(n^{-1})
\end{equation*}
where $F_b$ is defined at Equation \eqref{defbessel}.
\end{itemize}
\end{theorem*}

This theorem together with results of \cite{Le} lead to Proposition \ref{ultraspherical},
thus answering a question of M. Ledoux about the behavior of the suitably rescaled
largest eigenvalues.

This paper is organized as follows. Part \ref{prod-jac}
consists in explicit computations of densities.
Part \ref{sec-univ} gathers useful information about free probability
and Jacobi unitary ensembles, and establishes asymptotics for the
eigenvalues counting measures with free probabilistic tools.
Part \ref{sec-edge} 
provides asymptotics of suitably rescaled
kernels at the hard and soft edges, and inside the bulk of the spectrum.

{\it Acknowledgments. \it}
The author is currently a JSPS postdoctoral fellow at the university of Kyoto.
The results of section \ref{inside} of this paper were obtained during his PhD and
he acknowledges useful conversations with his advisor P. Biane, and also 
with T. Duquesne and J-F. Quint at an  early stage of the paper.
I also acknowledge stimulating discussions with M. Capitaine, M. Casalis,
Y. Doumerc and M. Ledoux about Jacobi unitary ensembles, and with A. Kuijlaars
about universality questions.

\section{Product of two random projections and Jacobi unitary ensembles.}\label{prod-jac}

Let $\mathbb{U}_n$ be the group of $n\times n$ complex unitary matrices,
and $\mu_n$ its normalized Haar measure. 
For $(\alpha ,\beta )\in\mathbb{R}^+$, consider
the probability distribution on the Hermitian matrices $\M{n}_{sa}$ given by
\begin{equation}\label{defjac}
(Z_n^{\alpha,\beta})^{-1}\det (1-M)^{\alpha}\det(M)^{\beta}1_{0\leq M\leq 1}dM
\end{equation}
where $Z_n^{\alpha,\beta}$ is some normalization constant.
This probability measure is
called {\it Jacobi unitary ensemble\it} of parameter $(n,\alpha,\beta )$
(see for example \cite{Fo}).

For $n,q$ positive integers and $a>0$, 
let $\tilde{W}(n,q,a)$ be the probability distribution on $\M{n\times q}$ whose density
is proportional to $e^{-a^{-1}\Tr (AA^*)}$.
Let $W(n,q,a)$ be the probability distribution on $\M{n}_{sa}$ of $WW^*$ where
$W\in\M{n\times q}$ has distribution $\tilde{W}(n,q,a)$. This probability measure is
called {\it Wishart ensemble\it} distribution and is proportional to
$\det (X)^{q-n}e^{-a^{-1}\Tr (XX^*)}dX$ whenever $q\geq n$.
We start with a classical lemma whose proof 
was explained to us by M. Casalis.

\begin{lemma}\label{lem2}
Let $X_{p},X_{p'}\in\M{n}$ be independent random matrices 
of distribution $W(n,p,1/n)$ and $W(n,p',1/n)$.
Then $X_{p}+X_{p'}$ is a Wishart matrix of parameter $(n,p+p',1/n)$.
Moreover, if $p+p'\geq n$, then almost surely, $X_{p}+X_{p'}$ is
invertible and we can define 
$$J= (X_{p}+X_{p'})^{-1/2} X_{p}(X_{p}+X_{p'})^{-1/2}$$
If $p,p'\geq n$, then the distribution of 
$J$ admits a density with
respect to the Lebesgue measure, and
it has the distribution of a Jacobi unitary ensemble of parameter $(n,p-n,p'-n)$.
\end{lemma}

\begin{proof}
Assume $p,p'\geq n$.
The random vector $(X_{p},X_{p'})$ has distribution
$$Ce^{-n\Tr (X+Y)}\det (X)^{p-n}\det (Y)^{p'-n}dXdY$$
By change of variable formula
together with the fact that the change of variable
$(X,Y)\rightarrow (X,X+Y)$ has Jacobian $1$,
$(X_p,X_p+X_{p'})$ has distribution
$$Ce^{-n\Tr (Y)}\det (Y-X)^{p'-n}\det (X)^{p-n}dXdY$$
The change of variable
$(X,S)\rightarrow (S^{-1/2}XS^{-1/2},S)$ on the cone of 
positive definite matrices is well defined, and
has Jacobian $\det (S)^d$. This implies that
$((X_p+X_{p'})^{-1/2}X_p(X_p+X_{p'})^{-1/2},X_p+X_{p'})$
has distribution
$$Ce^{-n\Tr (Y)}\det (Y)^{p-n}\det (1-X)^{p-n}\det (X)^{p'-n}dXdY$$
This proves that
$(X_p+X_{p'})^{-1/2}X_p(X_p+X_{p'})^{-1/2}$ 
has the distribution of a Jacobi unitary ensemble of parameter $(p-n,p'-n)$.
\end{proof}

\begin{theorem}\label{jacobi}
Let $X,X'\in\M{q_n}$ be independent Wishart matrices 
of parameter 
$(q_n,n-\tilde q_n,1/q_n)$ and $(q_n,\tilde q_n,1/q_n)$. 
Define $J$ as in Lemma \ref{lem2}, by
$$J=(X+X')^{-1/2}X(X+X')^{-1/2}$$
Let $\pi_n\in\M{n}$ be a constant orthogonal projection of rank $q_n$ and
$\tilde\pi_n\in\M{n}$ be a random uniformly distributed orthogonal projection of
rank $\tilde q_n$.
Then, under the unitary isomorphism
$\M{q_n}=\pi_n\M{n}\pi_n$, the following equality in distribution holds:
$$\pi_n\tilde\pi_n\pi_n\overset{\mathcal{L}}{=}J$$
In particular for $\tilde q_n\geq q_n$ and $q_n+\tilde q_n \leq n$,
$\pi_n\tilde{\pi}_n\pi_n$ has the distribution of a
 Jacobi unitary ensemble of parameter $(q_n,n-q_n-\tilde q_n, \tilde q_n -q_n)$ 
 (on $\M{q_n}$).
\end{theorem}

\begin{proof}
Let $\pi$ be a (deterministic) projection of rank $q_{n}$, $W$ and $\tilde{\pi}$ be independent
random matrices of $\M{n}$ having respective distributions $\tilde{W}(n,n,q_{n}^{-1})$ and
the invariant distribution on the selfadjoint  projectors of rank $\tilde{q}_{n}$.
Define $X_{1},X_{2}$ as
\begin{eqnarray}
X_{1 }& = & \pi W\tilde{\pi} W^{*}\pi \\
X_{2} & = & \pi  W(Id -\tilde{\pi}) W^{*}\pi
\end{eqnarray}
By construction, $X_{1}$ and $X_{2}$ are independent Wishart matrices in $\pi\M{n}\pi$.

Let $U$ be an unitary random variable such  that 
$$\pi W U^*=(X_{1}+X_{2})^{1/2}$$
This random variable can be chosen to depend measurably on $W$.
We have by definition
\begin{equation}
\label{11}
X_{1}=(X_{1}+X_{2})^{1/2}U\tilde{\pi}U^{*}(X_{1}+X_{2})^{1/2}
\end{equation}
therefore, from Lemma \ref{lem2},
$\pi U\tilde{\pi} U^{*}\pi$ has the distribution of a Jacobi unitary ensemble
of parameter $(q_{n},n-q_{n}-\tilde{q}_{n},\tilde{q}_{n}-q_{n})$.
Since $U$ is independent from $\tilde{\pi}$ (indeed, $U\in\sigma (W)$) and
$\tilde{\pi}$ is uniformly distributed, the distribution of $U\tilde{\pi} U^{*}$ is the same
as that of $\tilde{\pi}$.

Consequently, $\pi\tilde{\pi}\pi$ has also the distribution of a Jacobi unitary ensemble
of parameter $(q_{n},n-q_{n}-\tilde{q}_{n},\tilde{q}_{n}-q_{n})$.

\end{proof}

\begin{remark}
The hypothesis $q_n\leq \tilde q_n$ and $q_n+\tilde q_n\leq n$ is a necessary and sufficient condition for the distribution of $\pi_{n}\tilde{\pi_{n}}\pi_{n}$
 to admit a density with respect to the Haar measure of  $\pi_{n}\M{n}\pi_{n}$.
However, this case enables us to study the distribution of the eigenvalues set of
any product of the type
$\pi_n\tilde\pi_n\pi_n$ without any assumption on $q_n$ and $\tilde q_n$. Indeed,
 
1/ Assume $\tilde q_n+q_n\leq n$, but $\tilde q_n<q_n$. 
Since $\pi_n\tilde \pi_n\pi_n$ and $\tilde \pi_n\pi_n\tilde\pi_n$ are
unitarily conjugate to each other, the study of non-trivial eigenvalues 
of $\pi_n\tilde \pi_n\pi_n$ is equivalent to the study of those 
of $\tilde \pi_n\pi_n\tilde\pi_n$, and the latter is
a Jacobi unitary ensemble.

2/ Assume that $q_n+\tilde q_n >n$ and $\tilde q_n\geq q_n$. 
Then consider the conjugate random projector 
$\pi_n''=1-\tilde \pi_n$ of rank $q_n''=n-\tilde q_n$. 
The non-trivial eigenvalues of the ensemble
$\pi_n\pi_n''\pi_n$ are the image by the reflexion of 
center $1/2$ of the non-trivial eigenvalues of
$\pi_n\tilde \pi_n\pi_n$. One has $q_n''<q_n$ and 
$q_n+q_n''\leq n$ so we come back to case 1/

3/ Assume $q_n+\tilde q_n>n$ and $\tilde q_n<q_n$. 
Then with the notations of 2/, one has $q_n>q_n''$ and
$n<q_n+q_n''$. Therefore $\pi_n''\pi_n\pi_n''$ is in case 2/. 
Thus, we have showed how to handle any case.
\end{remark}

\section{Global asymptotics for Jacobi unitary ensembles.}\label{sec-univ}

\subsection{A reminder of free probability.}

We define a {\it non-commutative probability space \it} as
an algebra with unit endowed with a tracial 
state $\phi$. 
We denote such a space by $(A, \phi )$. An element of this space is called
a (non-commutative) random variable. 

Let $A_1, \cdots ,A_k$ be subalgebras of $A$ having the same unit as $A$.
They are said to be {\it free \it} iff
for all $a_i\in A_{j_i}$ ($i\in[1,k]$) 
such that $\phi(a_i)=0$, one has  
$$\phi(a_1\cdots a_l)=0$$
as soon as $j_1\neq j_2$, $j_2\neq j_3,\cdots ,j_{l-1}\neq j_l$.
Collections $S_{1},S_{2},\ldots $ of random variables are said to be 
{\it free \it} iff the 
unital subalgebras that they generate are free.

Let $(a_1,\cdots ,a_k)$ be a $k$ -tuple of random variables and let
$\mathbb{C}\langle X_1 , \cdots , X_k \rangle$ be the
free algebra of non commutative polynomials on $\mathbb{C}$ generated by
the $k$ indeterminates $X_1, \cdots ,X_k$. 
The {\it joint distribution\it} of the family $a_i$ is the linear form
$$\mu_{(a_1,\cdots ,a_k)} : 
\mathbb{C}\langle X_1, \cdots ,X_k \rangle
\rightarrow \mathbb{C}$$ 
defined in the obvious sense. 

Given a $k$ -tuple $(a_1,\cdots ,a_k)$ of free 
random variables and given each distribution $\mu_{a_i}$, the joint distribution
$\mu_{(a_1,\cdots ,a_k)}$ is uniquely determined by the
$\mu_{a_i}$'s.
A family $(a_1^n,\cdots ,a_k^n)_n$ of $k$ -tuples of random
variables {\it converges in distribution \it} towards $(a_1,\cdots ,a_k)$
iff for all $P\in \mathbb{C}\langle X_1, \cdots ,X_k \rangle$, 
$\mu_{(a_1^n,\cdots ,a_k^n)}(P)$ converges towards
$\mu_{(a_1,\cdots ,a_k)}(P)$ as $n\rightarrow\infty$. 
A sequence of families
$(a_1^n,\cdots ,a_k^n)_d$ is {\it asymptotically free \it} as $d\rightarrow\infty$
iff it converges in distribution towards a free random variable. 
Asymptotic freeness of sequence of collections of random variables is defined in
an analogous obvious sense.

The following result was contained in
 \cite{MR2000d:46080} and in  \cite{MR99f:81185} under slightly stronger
hypotheses. For a proof in full generality, see \cite{Co}, 
Proposition 2.3.3 p.52 or \cite{imrn}, Theorem 3.1.

\begin{theorem}\label{libre}
Let $U_{1},\cdots ,U_{k}, \cdots$ be a collection of independent 
Haar distributed random matrices of $\M{n}$ and $( W^n_i )_{i\in I}$ be a 
set of constant matrices of $\M{n}$ 
admitting a joint limit distribution for large $n$ with respect to the
state $n^{-1}\Tr$.
Then the family $( ( U_{1}, U_1^* ) ,\cdots , 
( U_{k}, U_k^* ), \cdots , (W_i))$ admits a limit distribution, and
is asymptotically free with respect to $E (n^{-1}\Tr )$.
\end{theorem}

\subsection{Free projectors.}

Let us fix real numbers 
$0\leq a \leq b \leq 1/2$, and let for all $n$, $\pi_{n}$ be a
self adjoint projector of $\M{n}$ of rank $q_{n}$ such that asymptotically
$q_{n}\sim \alpha n$ as $n\rightarrow\infty$.
Let $\pi_{n}'$ be a projector of rank $q_{n}'$ such that
$q_{n}'\sim \beta n$, and assume that it can be written under the form
$U\pi U^*$ such that $U$ is unitary Haar distributed
independent from $\pi_{n}$

It is a consequence of Theorem \ref{libre}, 
that $\pi_{n}$ and $\pi_{n}'$ are asymptotically free. 
Therefore $\pi_{n}\pi_{n}'\pi_{n}$ has an empirical eigenvalues
distribution converging towards $\mu_1\boxtimes\mu_2$, where $\mu_1$ is
the probability 
\begin{equation*}
(1-\alpha)\delta_0+\alpha\delta_1
\end{equation*}
and $\mu_2$ is
the probability 
\begin{equation*}
(1-\beta)\delta_0+\beta\delta_1
\end{equation*}
Let 
$$r_{\pm}=\alpha +\beta-2\alpha\beta\pm\sqrt{4\alpha\beta (1-\alpha)(1-\beta)}$$
By a standard $S$-transform argument (see \cite{MR94c:46133}, example
3.6.7), 

\begin{equation*}
\mu_1\boxtimes\mu_2=[1-\min (\alpha,\beta)]\delta_0+
[\max (\alpha +\beta-1,0)]\delta_1+
\frac{\sqrt{(r_+ -x)(x-r_-)}}{2\pi x(1-x)}1_{[r_-,r_+ ]}dx
\end{equation*}
By Theorem \ref{jacobi}, we recover a short proof of the following result:
\begin{proposition}[\cite{CC}, Corollary 7.2.]\label{propcc}
Let $X_n$ and $X_n'$ be independent complex 
Wishart matrices with respective distributions
$W(n,p_n,Id/n)$ and $W(n,p_n',Id/n)$, such that $p_n/n\sim \alpha\geq 1$ and
$p_n'/n\sim \beta\geq 1$, 
let $Z_n=(X_n+X_n')^{-1/2}X_n(X_n+X_n')^{-1/2}$. 
The expectation of the normalized eigenvalues counting measure tends
in moments towards
\begin{equation*}
\nu_{\alpha,\beta}(dx)=g(x)1_{[\lambda_-,\lambda_+]}dx
+\max (0,\alpha -1)\delta_0+\max (0,\beta -1)\delta_1
\end{equation*}
where
\begin{equation*}
\begin{split}
g(x)=\frac{\sqrt{(x-\lambda_-)(\lambda_+-x)}}{2\pi x(1-x)}\\
\lambda_{\pm}=\left(\sqrt{\frac{\alpha}{\alpha+\beta}(1-\frac{1}{\alpha+\beta})}
\pm\sqrt{\frac{1}{\alpha+\beta}(1-\frac{\alpha}{\alpha+\beta})}\right)^2
\end{split}
\end{equation*}
\end{proposition}

Indeed, by Theorem \ref{jacobi}, $Z_n$ has the same asymptotic distribution as
$\pi_{n}\pi_{n}'\pi_{n}$ where $\pi_{n},\pi_n'\in\M{p_n+p_n'}$ have
respective ranks $n$ and $p_n$, and the proposition follows by
a change of variables.  
Note that in \cite{CC} it is also proved  that
$Z_{n}$ is asymptotically free with $X_{n}+X_{n}'$.

\subsection{Jacobi polynomials and Jacobi kernel.}
In this section we gather technical results for the computation of asymptotics.

\subsubsection{Determinantal point process}\label{detpp}
Denote by $\Lambda_n^{a,b}=\{\lambda_1\geq\ldots \geq \lambda_n\}$
the random set of eigenvalues of a Jacobi unitary ensemble of parameter
$(n,a,b)$. 
Almost surely, this ensemble has cardinal $n$.

It is a so-called {\it determinantal point process\it}, i.e., there exists a kernel
$K_{n}^{a,b}$ which we will describe at section \ref{defker}, such that:
\begin{equation*}
\begin{split}
P(\Lambda_n^{a,b}\cap [x_1,x_1+dx_1]=1,\ldots ,\Lambda_n^{a,b}
\cap [x_1,x_1+dx_1]=1)=\\
dx_1\ldots dx_n\det (K_n^{a,b}(x_i,x_j))
\end{split}
\end{equation*}
We refer to \cite{MR2000g:47048}, and to \cite{MR92f:82002}
for a probabilistic interpretation.
Furthermore, $P (\lambda_1 \leq x )$
can be computed explicitely and its value is
(see Equation (5.42) p. 114 of \cite{MR2000g:47048})
\begin{equation}\label{deift}
\begin{split}
P (\lambda_1 \leq x )
=\det (I-K_n^{a,b})_{[x,\infty )}=\\
\sum_{j=0}^{\infty}\frac{(-1)^j}{j!}
\int_{[x,\infty ]}\ldots \int_{[x,\infty ]}
\begin{vmatrix}
K_n^{a,b}(x_1,x_1) & \cdots & K_n^{a,b}(x_1,x_j) \\
\vdots & & \vdots \\
K_n^{a,b}(x_j,x_1) & \cdots & K_n^{a,b}(x_j,x_j)
\end{vmatrix}
dx_1\ldots dx_j
\end{split}
\end{equation}

\subsubsection{The polynomials $P_{n}^{a,b}$ }\label{kernelk}
For $a,b \geq 0$, 
the Jacobi polynomials
$(P_n^{a,b})_{n\geq 0}$ form
a sequence of orthogonal polynomials with respect to the measure
$\frac{\Gamma (a+b+2)}{2^{a+b+1}\Gamma (a+1)\Gamma (b+1)}
1_{[-1,1]}w^{a,b}(x)dx$ where
\begin{equation*}
w^{a,b}(x)=(1-x)^{a}(1+x)^{b}
\end{equation*}
The normalization constant is such that (see \cite{MR51:8724}, Equation (4.3.4))
$$P_n^{a,b}(1)={n+a\choose n}$$

\begin{remark}
Observe that we choose to consider the weight
$(1-x)^{a}(1+x)^{b}$ instead of
$(1-x)^{a}x^{b}$ in order to respect the conventional notation
for Jacobi polynomials. The map $x\rightarrow 2x-1$ 
turns the Jacobi unitary ensemble $const \det (1-M)^{a}\det (M)^{b}dM$
obtained in Theorem
\ref{jacobi} into a ``classical'' Jacobi ensemble 
$const \det (1-M)^{a}\det (1+M)^{b}dM$. Therefore asymptotics of both ensembles
deduce from each other through elementary (affine~!) functional calculus.
\end{remark}
It will be useful to know that (see Szeg\"o \cite{MR51:8724})
\begin{equation}\label{contour}
P_n^{a,b}(x)(1-x)^{a}(1+x)^{b}
=\frac{1}{2i\pi}\int_{\Gamma}\frac{(1-t)^{a+n}(1+t)^{b+n}}{2^n(x-t)^n}
\frac{dt}{t-x}
\end{equation}
where $\Gamma$ is a closed $C^1$ curve with winding number $1$ around $x$ and
$0$ around $-1$ and $1$.
The formula (4.21.7) of \cite{MR51:8724} reads
\begin{equation}\label{j2}
(P_n^{a,b })'(x)=\frac{1}{2}(n+a+b+1)P_{n-1}^{a+1,b+1}(x)
\end{equation}
and formula (22.6.3) p781 of \cite{MR94b:00012} implies that
the function 
\begin{equation}\label{defgn}
g_n^{a,b}(x)=(1-x)^{(a+1)/2}(1+x)^{(b+1)/2}P_n^{a,b}(x)
\end{equation}
satisfies the second order differential equation
\begin{equation}\label{eqdif}
\frac{d^{2}}{dx^{2}}g_n^{a,b}=\chi_n^{a,b}g_n^{a,b}
\end{equation}
where
\begin{equation}\label{defeqdif}
\chi_n^{a,b}(x)=\frac{1-a^2}{4(1-x)^2}+\frac{1-b^2}{4(1+x)^2}
+\frac{2n(n+a+b+1)+(a +1)(b +1)}{2(1-x^2)}
\end{equation}
will be of fundamental use for our purposes.

\subsubsection{The kernel $K_{n}^{a,b}$}\label{defker}

Let $p_n^{a,b}$ be the corresponding orthonormal polynomials.
The kernel $K_{n}^{a,b}$ can be defined as the function
\begin{equation*}
K^{a,b}_n(x,y)=\sqrt{w^{a,b}(x)}\sqrt{w^{a,b}(y)}
\sum_{j=0}^{n-1}p^{a,b}_j(x)p^{a,b}_j(y)
\end{equation*}
which satisfies for $x\neq y$ by the Christoffel-Darboux formula
(see \cite{MR51:8724}, formula (4.5.2)):
\begin{equation}\label{j1}
K_n^{a,b}(x,y)=\gamma_n^{a,b}
\sqrt{w^{a,b}(x)}\sqrt{w^{a,b}(y)}
\frac{P_n^{a,b}(x)P_{n-1}^{a ,b}(y)-P_{n-1}^{a,b}(x)
P_n^{a,b}(y)}{x-y}
\end{equation}
Where
\begin{equation}\label{defGn}
\gamma_n^{a,b}=
\frac{2^{-a-b}}{2n+a+b}
\frac{\Gamma(n+1)\Gamma(n+a+b+1)}{\Gamma(n+a)\Gamma(n+b)}
\end{equation}
A Taylor expansion and Equation \eqref{j2} show that
\begin{equation}\label{one-point}
\begin{split}
K_n^{a,b}(x,x)=\gamma_n^{a,b}
w^{a,b}(x)
(P_{n-1}^{a,b}(x)(P_n^{a,b}(x))'-P_n^{a,b}(x)(P_{n-1}^{a,b}(x))')=\\
\frac{(n+a+b)}{2}[P_{n-1}^{a,b}(x)P_{n-1}^{a+1,b+1}(x)-
P_{n}^{a,b}(x)P_{n-2}^{a+1,b+1}(x)]\\
+\frac{1}{2}P_{n-1}^{a,b}(x)P_{n-1}^{a+1,b+1}(x)
\end{split}
\end{equation}

The function $x\rightarrow n^{-1}K_q^{a,b}(x,x)$ is the 
expectation of the normalized eigenvalues counting measure
of a Jacobi unitary ensemble of parameter $(a,b )$ on $\M{q}$ (see for example
\cite{MR92f:82002}, A.10). It is often called the 
{\it one point distribution function.\it}

\section{Local asymptotics and universality.}\label{sec-edge}

We make the following
\begin{assumption}\label{hypo}
Let $\alpha_n=a_n/n$ and $\beta_n=b_n/n$.
\begin{equation}
\lim_{n=\infty} \alpha_n =\alpha\in [0,+\infty) \, ,\,
\lim_{n=\infty} \beta_n =\beta\in [0,+\infty)
\end{equation}
\end{assumption}

Define 
\begin{equation}\label{defrsn}
\begin{split}
A_n=\frac{\alpha_n}{2+\alpha_n+\beta_n},\,\,\, B_{n}=\frac{\beta_n}{2+\alpha_n+\beta_n}\\
D_n=\sqrt{(1+A_n+B_n)(1-A_n-B_n)(1-A_n+B_n)(1+A_n-B_n)}\\
\tilde{r}_n=B_n^2-A_n^2-D_n,\,\,\, \tilde{s}_n=B_n^2-A_n^2+D_n
\end{split}
\end{equation}

In addition, let 
\begin{equation}\label{defrs}
A=\lim_{n}A_{n},B=\lim_{n}B_{n},
r=\lim_{n}\tilde{r}_{n}, s=\lim_{n}\tilde{s}_{n}
\end{equation}

Theorems  \ref{jacobi}, \ref{libre} and Proposition \ref{propcc}
imply 
\begin{corollary}
The probability measure
\begin{equation*}
n^{-1}K_n^{a_n,b_n}(x,x)dx
\end{equation*}
tends towards 
\begin{equation*}
1_{[r,s]}\frac{\sqrt{(x-r)(s-x)}}{\pi (1-A-B)(1-x^2)}dx
\end{equation*}
in the sense of moments.
\end{corollary}

This allows us to define 
\begin{equation}\label{deff}
\begin{split}
f_n(x)=\frac{\sqrt{(x-\tilde{r}_n)(\tilde{s}_n-x)}}{\pi (1-A_n-B_n)(1-x^2)}\\
f(x)=\frac{\sqrt{(x-r)(s-x)}}{\pi (1-A-B)(1-x^2)}
\end{split}
\end{equation}

We will see later (Theorem \ref{onepoint} and Proposition \ref{banach}) 
that the convergence of density functions 
actually holds uniformly on any compact set containing neither $r$ nor $s$.

\begin{remark}
The distribution $f_n(x)dx$ already appeared in the study of zeros and 
asymptotics of Jacobi polynomials (see \cite{MR80g:33021}).
Theorem \ref{jacobi} provides a simple 
explanation for the apparition of the same distribution in two a priori
very different places of mathematics.
\end{remark}

It is widely believed that ``reasonable'' unitary ensembles should have
universality properties for local spacing both inside the bulk of their
asymptotic spectra and at the edge up to some suitable renormalisation.
Amongst very recent results towards these conjectures, see the recent
work of Ledoux \cite{Le}. 

We settle this universality problem in the 
specific framework of Jacobi unitary ensemble satisfying hypothesis
\ref{hypo}.

Our approach is mainly based on the Christoffel-Darboux formula \eqref{j1}
and it only holds for non equal parameters in the kernel.
In order to settle this problem we will need the following reformulation
of the analytic maximum modulus principle:
\begin{lemma}\label{principe-maximum} 
Let $F_n : \mathbb{C}\times\mathbb{C}\rightarrow \mathbb{C}$ be a sequence
of holomorphic functions in both variables converging towards some
function $F$ uniformly on any compact subset of 
$\mathbb{C}^2-\{(x,x),x\in\mathbb{C}\}$. 
Then the limit $f$ extends by continuity to an holomorphic function
on $\mathbb{C}^2$, and the convergence holds on any compact subsets
of $\mathbb{C}^2$.
\end{lemma}

\subsection{Universality inside the spectrum.}\label{inside}

We first recall the following result of  \cite{MR92g:33011}.
Assume $a,b\in\mathbb{R}$, $\alpha,\beta\in\mathbb{R}_+$ and let

\begin{equation}
\Delta = [\alpha (x+1)+\beta (x-1)]^2-4(1+\alpha+\beta)(1-x^2)
\end{equation}

The nature of the asymptotics of $P_n^{\alpha n +a,\beta n+b}$ 
depends on the sign of $\Delta$. In the case $\Delta <0$,
for $x\in (-1,1)$, let $\rho,\theta,\gamma\in (-\pi ,\pi]$ be defined by

\begin{equation}\label{norme-tleft}
\frac{\alpha(x+1)+\beta(x-1)+i\sqrt{-\Delta}}{(1+\alpha+\beta)(1-x^2)}=
 2[(1+\alpha+\beta)(1-x^2)]^{-1/2}e^{i\rho} 
 \end{equation}
 
 \begin{equation}\label{norme-xi}
\frac{(\alpha+\beta+2)x-(3\alpha+\beta+2)-i\sqrt{-\Delta}}{2(x-1)(\alpha+\beta+1)}=
\left[\frac{2(\alpha+1)}{(1-x)(\alpha+\beta+1)}\right]^{1/2}e^{i\theta}
\end{equation}

\begin{equation}\label{norme-eta}
\frac{(\alpha+\beta+2)x+(3\alpha+\beta+2)-i\sqrt{-\Delta}}{2(x+1)(\alpha+\beta+1)}=
\left[\frac{2(\beta+1)}{(1+x)(\alpha+\beta+1)}\right]^{1/2}e^{i\gamma}
\end{equation}

The result is 

\begin{equation}\label{chis}
\begin{split}
P_n^{\alpha n+a,\beta n+b}=\left(\frac{4\sqrt{-\Delta}}{\pi n}\right)^{1/2}
\left[\frac{(1-x)(\alpha+\beta+1)}{2(\alpha+1)}\right]^{n(-\alpha-1)/2-a/2-1/4}\\
\left[\frac{(1+x)(\alpha+\beta+1)}{2(\beta+1)}\right]^{n(-\beta-1)/2-b/2-1/4}
\left[\frac{(1-x^2)(\alpha+\beta+1)}{4}\right]^{n/2+1/4}\\
(\cos [((1+\alpha)n+a +1/2)\theta/2+(-(1+\beta)n+b +1/2)\gamma/2\\
-(2n+1)
\rho/4-\pi/4]+(O(1/n)))
\end{split}
\end{equation}

It is a consequence of \cite{MR92j:33019,MR2001i:30040} 
(see also \cite{Co}, Lemma 4.3.2 pp. 120-121 )
that this estimate is uniform in compact subsets
$K$ of $\mathbb{R}^+\times\mathbb{R}^+\times \mathbb{R}\times \mathbb{R}
\times (-1,1)$ such that for any element in $K$, the associated $\Delta$ is 
negative.

Let $\gamma_n=\gamma_n^{a_n,b_n}$ where $\gamma_n^{a,b}$ was defined 
at Equation \eqref{defGn}.
The following is a straightforward application from Stirling's asymptotic formula.
\begin{lemma}\label{stirling}
As $n\to\infty$,
\begin{equation*}
\gamma_n
=n\frac{2^{-a_n-b_n}(1+\alpha_n+\beta_n)^{n+a_n+b_n+1/2}}{(1+\alpha_n)^{n+a_n-1/2}(1+\beta_n)^{n+b_n-1/2}(2+\alpha_n+\beta_n)} (1+O(n^{-1}))
\end{equation*}
\end{lemma}

\begin{theorem}\label{onepoint}
The following holds true:
\begin{equation*}
n^{-1}K_n^{a_n,b_n}(x,x)= f(x)(1+O(n^{-1}))
\end{equation*}
where $O(n^{-1})$ holds uniformly on compact subsets of 
$(r,s)$.
(see the definition of $r,s$ at Equation \eqref{defrsn}),
$f$ was defined at Equation \eqref{deff}.
\end{theorem}

\begin{proof}
Recall that
\begin{equation*}
\begin{split}
n^{-1}K_n^{a_n,b_n}(x,x)=\\
(1-x)^{a_n}(1+x)^{b_n}\gamma_n
(P_{n-1}^{a_n,b_n}(x)(P_n^{a_n,b_n}(x))'-P_n^{a_n,b_n}(x)(P_n^{a_n,b_n}(x))')
\end{split}
\end{equation*}
By Lemma \ref{stirling}, the right hand side is equivalent to
\begin{equation}\label{stirl}
\begin{split}
(1-x)^{a_n}(1+x)^{b_n}\frac{2^{-a_n-b_n}}{2+\alpha+\beta}
\frac{(1+\alpha+\beta)^{n+a_n+b_n+1/2}}{(1+\alpha)^{n+a_n-1/2}
(1+\beta)^{n+b_n-1/2}}\\
\frac{a_n+b_n+n}{2}
(P_{n-1}^{a_n,b_n}(x)P_{n-1}^{a_n+1,b_n+1}(x)-P_n^{a_n,b_n}(x)P_{n-2}^{a_n+1,b_n+1}(x))
\end{split}
\end{equation}
Plugging into Equation \eqref{stirl} the asymptotics of Formula \eqref{chis}
yields 
\begin{equation*}
\begin{split}
n^{-1}K_n^{a_n,b_n}(x,x)=
\frac{(\alpha_{n}+\beta_{n}+1)2\sqrt{-\Delta}}{\pi \sqrt{1-x^2}(2+\alpha_{n}+\beta_{n})
(1+\alpha_{n})(1+\beta_{n})}\times \\
(\cos (\frac{\theta_{n}-\gamma_{n}}{2})[1-2\cos \rho_{n}]-\sin (\frac{\theta_{n}-\gamma_{n}}{2})\sin\rho_{n} )
(1+O(n^{-1}))
\end{split}
\end{equation*}

where $\theta_{n},\rho_{n},\gamma_{n}$ are defined in
Equations \eqref{norme-tleft},\eqref{norme-xi},\eqref{norme-eta}
by replacing $\alpha$ by $\alpha_{n}$ and $\beta$ by $\beta_{n}$.
Direct computation shows that
$$(\cos (\frac{\theta-\gamma}{2})[1-2\cos \rho]
-\sin (\frac{\theta-\gamma}{2})\sin\rho )=
\frac{(1+\alpha)(1+\beta)}{2(\alpha+\beta+1)\sqrt{1-x^2}} $$
therefore we obtain
\begin{equation*}
n^{-1}K_n^{a_n,b_n}(x,x)= 
\frac{\sqrt{(x-\tilde{r}_n)(\tilde{s}_n-x)}}{\pi (1-A_n-B_n)(1-x^2)}+0(n^{-1})
\end{equation*}
and this completes the proof.
\end{proof}

Using the same method (involving cumbersome calculations with
Formula \eqref{chis}) to treat the universality in the bulk of
the spectrum in the same fashion, one finds:

\begin{theorem}\label{univ-local}
For $u,v\in (0,\infty)$, we have as $n\rightarrow\infty$,
\begin{equation*}
\frac{1}{nf_n(x)}K_n^{a_n,b_n}(x+\frac{u}{nf_n(x)},x+\frac{v}{nf_n(x)})=
\frac{\sin \pi (u-v)}{\pi(u-v)}+0(n^{-1})
\end{equation*}
This limit is uniform for $x$ in any compact  subset of 
$(r,s)$ and for $u,v$ in compact subsets of $\mathbb{R}$. 
\end{theorem}

Here, we present a proof which is more intrinsic and instructive for
asymptotics at the edges, making use of the differential equation
\eqref{eqdif}. However, with this proof we do not obtain the optimal
error term.

\begin{proof}[Proof with error term $O(n^{-1+\varepsilon})$]
It is a consequence of Equation \eqref{eqdif}
that for any $x$, the functions
\begin{equation*}
p^1_{n,x} : u\rightarrow
(1-x-\frac{\pi u}{nf_n(x)})^{(a_n+1)/2}(1+x+\frac{\pi u}{nf_n(x)})^{(b_n+1)/2}
P_n^{a_n,b_n}(x+\frac{\pi u}{nf_n(x)})
\end{equation*}
and
\begin{equation*}
p^2_{n,x} : u\rightarrow
(1-x-\frac{\pi u}{nf_n(x)})^{(a_n+1)/2}(1+x+\frac{\pi u}{nf_n(x)})^{(b_n+1)/2}
P_{n-1}^{a_n,b_n}(x+\frac{\pi u}{nf_n(x)})
\end{equation*}
satisfy the differential equation
\begin{equation}\label{difloc}
p^{*''}_{n,x}+(1+O(n^{-1}))p_{n,x}^*=0
\end{equation}
where $*\in \{1,2\}$ and
the term $o(n^{-1})$ has to be understood as $n\rightarrow\infty$
uniformly on any compact set of couples $(x,u)$ such 
that the constant $\Delta$ determined
by $x$ is strictly negative. 

Consider the analytic function
$u\rightarrow p^1_{n,x}(u)p^{2}_{n,x}(v)-p^2_{n,x}(u)p^{1}_{n,x}(v)$
Its value in $v$ is zero and its first order derivative $c_n$ is given by
Theorem \ref{onepoint}. Using Equation \eqref{difloc}, one can apply
a recursive approximation argument on the derivatives of this 
function in $v$ (see Lemma \ref{bound} for detail) to see that
for any $\varepsilon >0$, 
$$c_n^{-1}(p^1_{n,x}(u)p^{2}_{n,x}(v)-p^2_{n,x}(u)p^{1}_{n,x}(v))
=\sin (\pi (u-v))/\pi+O(n^{-1+\varepsilon })$$
uniformly on compact subsets of $\mathbb{C}$. This result can be checked
to hold uniformly for $v$ in compact sets of $\mathbb{C}$.
Therefore 
\begin{equation*}
\frac{u-v}{n}K_n^{a_n,b_n}(x+\frac{u}{nf_n(x)},x+\frac{v}{nf_n(x)})\sim 
f_n(x)\sin (\pi (u-v))/\pi
\end{equation*}
This implies the theorem on compact subsets of $\mathbb{C}^2$ such that
$u\neq v$. The general result on arbitrary compact sets of $\mathbb{C}^2$ follows
by Lemma \ref{principe-maximum}.
\end{proof}

\begin{proposition}[\cite{MR92g:33011}]\label{asymp}
Assume $\Delta > 0$.
Let 
\begin{equation}\label{deft}
t_{\pm}=\frac{\beta(x-1)+\alpha(1+x)\pm\sqrt{\Delta}}{(\alpha+\beta+1)(1-x^2)}
[1+\xi_{\pm}]^{-\alpha-1}[1+\eta_{\pm}]^{-\beta-1}
\end{equation}
where
\begin{equation}\label{defxi}
\xi_{\pm}=\frac{\beta(x-1)+\alpha(1+x)\pm\sqrt{\Delta}}{-2(\alpha+\beta+1)(x-1)}
\end{equation}
and
\begin{equation}\label{defeta}
\eta_{\pm}=\frac{(x-1)\xi_{\pm}}{x+1}
\end{equation}
Let $t_0=t_+$ if $|t_-|>|t_+|$ and $t_0=t_-$ if $|t_-|<|t_+|$.
Then there exists a non zero real number $C_0$ such that
\begin{equation}
|P_n^{\alpha n,\beta n}(x)|\leq\frac{C_0}{\sqrt{n}}t_0^{-n}
\end{equation}
\end{proposition}
Furthermore, it is a consequence of \cite{MR2001i:30040} that
this estimate is uniform on compact subsets of $\{x,\Delta >0\}$.

\begin{lemma}\label{ineq}
Let $\Delta \geq 0$
\begin{itemize}
\item
\begin{equation*}
\min ([1+\xi_{-}]^{2},[1+\xi_{+}]^{2})\leq\frac{2(\alpha+1)}{(1-x)(\alpha+\beta+1)}
\leq \max ([1+\xi_{-}]^{2},[1+\xi_{+}]^{2})
\end{equation*}
with equality iff $\Delta =0$
\item
\begin{equation}
\min ([1+\eta_{+}]^{2},[1+\eta_{-}]^{2})\leq
\frac{2(\beta+1)}{(1+x)(\alpha+\beta+1)}
\leq\max ([1+\eta_{+}]^{2},[1+\eta_{-}]^{2})
\end{equation}
with equality iff $\Delta =0$
\item
\begin{equation*}
\begin{split}
\min [
\left(\frac{\alpha(x+1)+\beta(x-1)+\sqrt{\Delta}}{(1+\alpha+\beta)(1-x^2)}\right)^2 ,
\left(\frac{\alpha(x+1)+\beta(x-1)-\sqrt{\Delta}}{(1+\alpha+\beta)(1-x^2)}\right)^2]\\
\leq 4[(1+\alpha+\beta)(1-x^2)]^{-1}\\
\leq
\max [
\left(\frac{\alpha(x+1)+\beta(x-1)+\sqrt{\Delta}}{(1+\alpha+\beta)(1-x^2)}\right)^2 ,
\left(\frac{\alpha(x+1)+\beta(x-1)-\sqrt{\Delta}}{(1+\alpha+\beta)(1-x^2)}\right)^2]
\end{split}
\end{equation*}
with equality iff $\Delta=0$
\end{itemize}
\end{lemma}

\begin{proof}[Proof (sketch)]
This is a consequence of \cite{MR92g:33011}. 
For the first point, it is a consequence of
Equations \eqref{norme-xi} and \eqref{defxi} and the complex triangular
inequality.
For the second (resp. third) inequality, make use of 
Equation \eqref{norme-eta} and \eqref{defeta}  (resp. \eqref{norme-tleft}).
\end{proof}

\begin{proposition}\label{norme}
Let $a_n$ and $b_n$ be sequences satisfying Assumption \ref{hypo},
and $\varepsilon\in (0,1-s)$. 

There exists a constant $C_1\in (0,1)$ depending on $\varepsilon$, $\alpha,\beta$ 
such that for $n$ large enough, for all $x\in [s+\varepsilon ,1]$, 
$K_n^{a_n,b_n}(x,x)\leq C_1^n$. 

As a consequence, almost surely, there is no eigenvalue in $[s+\varepsilon ,1]$ for
$n$ large enough in the sequence of Jacobi unitary ensembles of parameter
$(n,a_n,b_n)$.
\end{proposition}

\begin{proof}
For the first point, 
Proposition \ref{asymp} together with Lemma \ref{ineq} and 
Equation \eqref{stirl} show that according to the definition of the
one point correlation function of Equation \eqref{one-point}, one has
$|K_n^{a_n,b_n}(x,x)|\leq P(n)C_0^n$ where $P$ is some polynomial.
Therefore any $C_{1}\in (C_{0},1)$ satisfies the announced property.

For the second point, 
observe that the summand of Equation \eqref{deift} 
satisfies
$$
\begin{vmatrix}
K_n^{a,b}(x_1,x_1) & \cdots & K_n^{a,b}(x_1,x_j) \\
\vdots & & \vdots \\
K_n^{a,b}(x_j,x_1) & \cdots & K_n^{a,b}(x_j,x_j)
\end{vmatrix}
\leq n^nC_1^{n^2}
$$
therefore
$$P(\lambda_{1}\leq s+\varepsilon )\geq 2-\exp (n^nC_1^{n^2})$$ 
for $n$ large enough.
The result follows by the Borel-Cantelli Lemma
together with the fact that the series
$(P(\lambda_{1}\geq s+\varepsilon ))_{n}$ has finite sum.
\end{proof}

\begin{remark}
Proposition \ref{norme} show that the result of
Theorem \ref{onepoint} extends on any compact set 
remaining at positive distance from both $\tilde{r}_n$ and
$\tilde{s}_n$
(defined at Equation \eqref{defrsn}). In the case both $\tilde{r}_n$
and $\tilde{s}_n$ stay at a positive distance froom $1,-1$,
A. Kuijlaars asked
whether the uniformity holds on any compact subset of $(-1,1)$.
We have not been able to answer this question if the compact set
$K$ contains the transition points $r,s$. It seems that 
Riemann-Hilbert techniques (see \cite{math-ph/0204006} for a 
bibliography) could give more insight.
\end{remark}

Proposition \ref{norme} has an interesting consequence in 
terms of geometry of Hilbert spaces:

\begin{proposition}\label{banach}
Let $\Omega_{n,q}$ be the set of subspaces of $\mathbb{C}^n$ of dimension
$q$ together with its uniform probability measure $P_{q}$.
Assume that there exists $\eta >0$ such that for all $n$,
$\alpha_n+\beta_n<1-\eta$.

It is possible to find an angle $\theta\in [0,\pi /2)$ satisfying,
for any $\varepsilon >0$, the existence of a constant $c>0$
such that for all $n,q,q'$ satisfying
$q\leq \alpha n , q'\leq \beta n$,
there exists a subset $F$ of $\Omega_{n,q}\times\Omega_{n,q'}$ of
measure larger than $1-e^{-cn}$ such that
\begin{equation*}
\begin{split}
\forall (V_1,V_2)\in F, \forall x_1\in V_1-\{0\}, x_2\in V_2-\{0\},
angle (x_1,x_2)\in [\theta-\varepsilon ,\pi /2 ]
\end{split}
\end{equation*}
The value of $\theta$ is given by
$$\cos^2 \theta =s $$ 
where $s$ was defined in Equation \eqref{defrs}
(observe that the assumption on $\eta$ ensures the existence
of $\theta\in [0,\pi /2)$).
\end{proposition}

\begin{proof}
This is a consequence of Proposition \ref{norme}
together with the fact that $||\pi \pi '||^2=||\pi\pi '\pi ||$.
\end{proof}

\subsection{At the soft edge.}

Instead of assumption \ref{hypo},
we shall 
work under the following slightly different assumption.

\begin{assumption}\label{soft}
Let $\alpha_n=a_n/n, \beta_n=b_n/n$. Then
$$\liminf_n \alpha_n >0,\limsup_n \alpha_n<\infty$$
and
$$\liminf_n \beta_n \geq 0,\limsup_n \beta_n<\infty$$
\end{assumption}

The clause $\liminf_n \alpha_n >0$ ensures that the behavior of the kernel
in the neighborhood of $\tilde{s}_n$ will be of ``soft edge'' type.
In assumption \ref{hypo} it was rather natural to assume the $\liminf$ and
$\limsup$'s are actual limits. It was not necessary but 
allowed lighter notations (in particular $\alpha_n$ and $\beta_n$ could
be replaced by their limits in the formulas for asymptotics). 
At the edge of the spectrum, assuming that
the $\liminf$ and $\limsup$'s are actual limits do not simplify the notation
(for example replacing $\tilde{s}_n$ by its possible limit) because 
unlike for the asymptotics inside the bulk of the spectrum, a control on the
speed of convergence is also needed.

\subsubsection{Preliminary estimate.}

\begin{proposition}[Chen-Ismail, \cite{MR92g:33011}]\label{softestimate}
The following holds true
\begin{equation*}
\begin{split}
P_n^{a_n,b_n}(\tilde{s}_n)=(1+o(n^{-1/3}))
\frac{2^{1/3}\sqrt{(\alpha_n+1)(\beta_n+1)}}{9^{1/3}n^{1/3}\Gamma (2/3)
(\alpha_n+\beta_n+1)^{2/3}}\\
\left(\frac{2(\alpha_n +1)}{(1-\tilde{s}_n)(\alpha_n+\beta_n+1)}\right)^{\alpha_n /3}
\left(\frac{2(\beta_n +1)}{(1+\tilde{s}_n)(\alpha_n+\beta_n+1)}\right)^{\beta_{n} /3}\\
(1-\tilde{s}_n)^{-a_n/2-1/6}(1+\tilde{s}_n)^{-b_n/2-1/6}
\left(\frac{(1+\alpha_n)^{n+a_n}
(1+\beta_n)^{n+b_n}}{(1+\alpha_n+\beta_n)^{n+a_n+b_n}} \right)^{1/2}
\end{split}
\end{equation*}
\end{proposition}

\begin{remark}\label{steepest}
\begin{itemize}
\item
This result is contained in \cite{MR92g:33011} with different notations.
The proof of \cite{MR92g:33011} holds only for arithmetic sequences
$a_n,b_n$, but uniformity can be easily derived from the paper using
Bessel inequality (see \cite{Co}, Lemma 4.3.2 pp. 120-121 ).
Alternatively, it is possible to choose the approach of \cite{MR2001i:30040}
making use of Integral \eqref{contour}. 
\item
In this section we only handle the case of $\tilde{s}_n$ but the obvious 
counterparts of Assumption \ref{soft} at the edge $\tilde{r}_n$ 
also holds true.
\end{itemize}
\end{remark}

The {\it Airy Equation\it} 
$f''=xf$, has a conventional basis for its solutions,
denoted by $(Ai,Bi)$, where
\begin{equation*}
\pi Ai (x)=\lim_{u\rightarrow +\infty}\int_0^{u}\cos (t^3/3+xt)dt
\end{equation*}
and
\begin{equation*}
\pi Bi (x)=\lim_{u\rightarrow +\infty}\int_0^{u}[\exp (-t^3/3+xt)+\sin (t^3/3+xt)]dt
\end{equation*}
On $[0,\infty )$, $Ai$ is positive and tends towards zero as
$x\rightarrow\infty$, whereas $Bi$ is positive and increases towards infinity.
Besides, $Ai(0)=9^{-1/3}\Gamma(2/3)^{-1}$.
The {\it Airy Kernel \it} is defined as 
\begin{equation}\label{defairy}
Ai(u,v)=\frac{Ai(u)Ai'(v)-Ai(v)Ai'(u)}{u-v}
\end{equation}

\subsubsection{Computation of the kernel.}

Let $s_n$ be the largest zero of $\chi_n^{a_n,b_n}$ in $(-1,1)$
(where $\chi_n^{a_n,b_n}$ was defined in Equation \eqref{defeqdif})
and $h_n$ be the real number such that 
$-h_n^3$ is the derivative of $\chi_n^{a_n,b_n}$ at $s_n$.
The sequences $s_n$ and $h_n$ actually depend on $n,a_n,b_n$.
One checks that as $n\to\infty$,
\begin{equation}\label{ashn}
h_n =n^{2/3}\left(
\frac {\sqrt { \left( 1+\alpha_n \right)  \left( 1+\beta_n \right)  
\left( 1+\alpha_n+\beta_n \right) }}{2 \left( 1- s_n^2 \right) ^{2}}\right)^{1/3}
(1+0(n^{-1/3}))
\end{equation}
Let 
$$\phi_n(x)=g_n^{a_n,b_n}(s_n+xh_n^{-1})$$ 
where $g_{n}^{a_{n},b_{n}}$ was defined at Equation \eqref{defgn},
and
$$\tilde{\phi}_n(x)=g_{n-1}^{a_n,b_n}(s_n+xh_n^{-1})$$
For any real number $R$, the function $\phi_n$ is defined on the
interval $[R,n^{2/3}(1-s_n)]$ for $n$ large enough. Furthermore,
its value is zero in $h_n(1-s_n)$. 
Let 
\begin{equation}\label{tildechi}
\tilde{\chi}_n(x)=\chi^{a_n,b_n}(s_n+xh_n^{-1})h_n^{-2}
\end{equation}
The function $\tilde{\chi}_n$ is positive on the interval
$[0,h_n(1-s_n)]$ and tends uniformly on compact sets towards the
identity function. By construction, its value in zero is zero.
The function $\phi_n$ satisfies the differential equation
\begin{equation}\label{equadif}
\phi_n"=\tilde{\chi}_n\phi_n
\end{equation}

\begin{lemma}\label{decrease}
The functions $\phi_n$, $\tilde{\phi}_n$ are
decreasing and positive on $[0,h_n(1-s_n)]$.
\end{lemma}

\begin{proof}
We prove only the result for $\phi_n$, the proof for
$\tilde{\phi}_n$ being exactly the same.
The fact that $\phi_n(0)$ is positive is a consequence of Proposition
\ref{softestimate}.

Let us first show that for $n$ large enough, $\phi_n$ is decreasing
for $x\geq 0$.
Assume that $\phi_{n}$ has a zero inside
$(0,h_n(1-s_n))$. This implies that there exists $x$ such that 
$\phi_{n}(x)<0$. Therefore there exists $x$ such that
both
$\phi_{n}(x)<0$ and $f_{n}(x)'<0$.
Since $\phi_n$ satisfies the differential equation
\eqref{equadif} and since $\tilde{\chi}_n$ is strictly positive on
the interval $[0,h_n(1-s_n))$, this contradicts the fact that
$\phi_{n}(h_n(1-s_n))=0$.
Therefore $\phi_{n}$ has no zero inside $(0,h_n(1-s_n))$.
This implies that for $n$ large enough, $\phi_{n}$
is positive on $(0,h_n(1-s_n))$. 
If one assumes that it is strictly increasing at some place, this would again
contradict the fact that
$\phi_{n}(h_n(1-s_n))=0$ by Rolle's theorem.
Therefore we have proved that
for $n$ large enough, $\phi_{n}$ is decreasing
for $x\geq 0$.
\end{proof}

\begin{lemma}\label{bound}
$1/800<-\frac{\phi_n(0)'}{\phi_n(0)}<2$
\end{lemma}

\begin{proof}
By  Lemma \ref{decrease}, the function $\phi_n'$ is negative on $(0,h_n(1-s_n))$ and increasing,
therefore $|\phi_n'(1)|\leq \phi_n(0)$ (or else integrating $\phi_n'$ over the interval $[0,1]$ would
contradict $\phi_n$ ranging in $[0,\phi_n(0)])$.
Integrating $\phi_n"$ over $[0,1]$ and using the differential equation therefore shows
that $|\phi_n'(0)|\leq 2\phi_n(0)$.

For the other inequality, observe first that $\phi_n(1/10)\geq \phi_n(0)/2$. Indeed, if this
were not the case, by positivity assumption, 
there would be $t\in [0,1/10], -\phi_n'(t)\geq 5\phi_n(0)$ and some
$t'\in [1/10,6/10],-\phi_n'(t')\leq \phi_n$, which would result in the existence of
some $t''\in [0,6/10], \phi_n''(t'')\geq 4\phi_n(0)$. This contradicts the differential
equation because $\phi_n$ is decreasing.

On the interval $[1/20,1/10]$, 
$\phi_n''\geq \phi_n(1/10)/20$ thus by the preceding inequality, $\phi_n''\geq \phi_n(0)/40$.
Integrating, this yields $-\phi_n'(1/20)\geq -\phi_n'(1/10)+\phi_n(0)/800\geq \phi_n(0)/800$.
\end{proof}

\begin{lemma}\label{approxt}
For any $\varepsilon >0$, one has
$\phi_n(x)/\phi_n(0)-Ai(x)/Ai(0)=0(n^{-2/3+\varepsilon})$
uniformly on compact subsets of $\mathbb{C}$.
\end{lemma}

\begin{proof}
In this proof, denote the power series expansion of $\phi_n$ and $Ai$
by
$$\phi_n(x)=\sum_k\phi_n[k]x^k, Ai(x)=\sum_kAi[k]x^k$$
An application of Cauchy integral formula to Equation \eqref{contour}
on a circular contour of center $\tilde{s}_n$
and diameter $d(\tilde{s}_n,\{-1,1\})/2$,
and the integral triangle inequality shows that
that there exists a constant $C_1 >0$ such that
$$|\phi_n[k]|\leq C_1^kn^{-2/3(k-1)}$$
Therefore, 
for any $\eta >0, \eta'\in (0,2/3)$, there exists a constant $C_2>0$
$x\in [-n^{\eta'},n^{\eta '}]$
\begin{equation}\label{reste}
\sum_{k\geq n^{\eta}}\phi_n[k]x^k\leq C_2n^{-(2/3-\eta ')n^{\eta}}n^{2/3}
\end{equation}
Denote the power series expansion of $\tilde{\chi}_n(x)$ by
$\tilde{\chi}_n(x)=\sum_{k\geq 0}b_kx^k$,
where, by assumption $b_0=0, b_1=1$.
There exists a constant $C_3$ such that for all $k\geq 2$,
$$|b_k|\leq C_3^kn^{-2(k-1)/3}$$
According to the differential equation structure of Equation \eqref{equadif},
the following recursive equation is satisfied:
$$n(n-1)a_n=\sum_{i=1}^{n-2}b_ia_{n-2-i}$$
Thanks to this equation and with Lemma \ref{bound}, 
one can show that there exist constants
$\eta >0$ and $C_4>0$ such that for $n$ large enough, 
and for any $k\leq n^{\eta}$,
\begin{equation}\label{boundairy}
\begin{split}
|\phi_n[3k]/\phi_n[0]-Ai[3k]/Ai(0)|\leq C_4n^{-2/3}Ai[3k]\\
|\phi_n[3k+1]/\phi_n[1]-Ai[3k+1]/Ai(0)|\leq C_4n^{-2/3}Ai[3k+1]\\
|\phi_n[3k+2]|\leq C_4n^{-2/3}Ai[3k]
\end{split}
\end{equation}
Fix $\eta '\in (0,2/3)$.
Equations \eqref{boundairy} show that there exists a constant $C_5>0$
such that
\begin{equation*}
\begin{split}
\sum_{k=0}^{n^{\eta}} \phi_n[3k]n^{k\eta '}/\phi_n[0]=
\sum_{k=0}^{n^{\eta}} Ai[3k]n^{k\eta '}/Ai(0) (1+0(n^{-2/3}))\\
\sum_{k=0}^{n^{\eta}} \phi_n[3k+1]n^{k\eta '}/\phi_n[0]=
\sum_{k=0}^{n^{\eta}} Ai[3k+1]n^{k\eta '}/Ai(0) (1+0(n^{-2/3}))\\
\sum_{k=0}^{n^{\eta}} \phi_n[3k+2]n^{k\eta '}/\phi_n[0]\leq C_5n^{-2/3}
\sum_{k=0}^{n^{\eta}} Ai[3k]n^{k\eta '}/Ai(0)
\end{split}
\end{equation*}
Besides,  $\sum_{k=0}^{n^{\eta}} Ai[3k]n^{k\eta '}/Ai(0)$
and $\sum_{k=0}^{n^{\eta}} Ai[3k+1]n^{k\eta '}/Ai(0)$
grow quicker than any polynomial as $n\rightarrow\infty$
because all summands have the same sign.
This together with the remainder estimate \eqref{reste}
and the fact, by
Lemma \ref{decrease}, that one has $\phi_n(n^{\eta '})\in [0,\phi_n(0)]$,
imply that for any $\varepsilon >0$,
$$|\phi_n[1]/\phi_n(0)-Ai'(0)/Ai(0)|=0(n^{-2/3+\varepsilon})$$
An application of 
Inequalities \eqref{reste} and \eqref{boundairy} concludes the proof.
\end{proof}

\begin{lemma}
Let $c_n=2(2n+a_n+b_n+2)/h_n^2$. Then
$$\tilde{\phi}_n(x)/\tilde{\phi}_n(0)-Ai(x-c_n)/Ai(-c_n)=0(n^{-2/3+\varepsilon})$$
uniformly on compact subsets of $\mathbb{C}$.
\end{lemma}

\begin{proof}
Observe that $x\rightarrow\tilde{\phi}_n(x)$
satisfies the differential equation
$$\tilde{\phi}_n(s_n'+\frac{x}{h_n})"=
\tilde{\phi}_n(s_n'+\frac{x}{h_n})\chi_{n-1}^{a_n,b_n}(s_n'+\frac{x}{h_h})h_n^{-2}$$ 
where $s_n'$ is the largest zero of $\chi_{n-1}^{a_n,b_n}$ on $(-1,1)$.
Besides, the function
$$x\rightarrow\chi_{n-1}^{a_n,b_n}(s_n'+x/h_h)/h_n^2 $$
has the same properties as $\tilde{\chi}_n$, namely it
tends uniformly on any compact set towards identity function
and is positive on $\mathbb{R}_+$.

An application of Taylor approximation formula 
shows that 
$$s_n'=s_n+c_n+O(n^{-2/3})$$
therefore we are exactly in the hypotheses of
Lemma \ref{approxt} and the proof follows in the same way.
\end{proof}

\begin{lemma}\label{diff}
Uniformly on compact subsets of $\mathbb{C}$, for any
$\varepsilon >0$, we have
$\tilde{\phi}_n(x)-\phi_n(x)=c_n\phi_n'(x)+0(n^{-2/3+\varepsilon})$.
\end{lemma}

\begin{proof}
First observe that
$\tilde{\phi}_n(x)=\phi_n(x-c_n)+0(n^{-2/3})$ by the previous Lemma.
Then it is standard that
$$\phi_n(x-c_n)-\phi_n(x)=c_n\phi_n'(x)+0(n^{-2/3})$$
by standard power series analysis.
\end{proof}

Therefore we end up with
$$\bar{f}_n(x)\phi_n(y)-\bar{f}_n(y)\phi_n(x)=c_n (\phi_n(x)\phi_n'(y)-\phi_n(y)\phi_n'(x))+0(n^{-2/3+\varepsilon})$$
By Formula \eqref{j1} we have, for $x\neq y$.
$$\frac{(x-y)}{h_n}K_n^{a_n,b_n}(s_n+\frac{x}{h_n},s_n+\frac{y}{h_n})
=\frac{\gamma_n}{1-s_n^2}(\bar{f}_n(x)\phi_n(y)-\bar{f}_n(y)\phi_n(x))$$
which by Lemma \ref{diff} is
$$\frac{\gamma_n}{1-s_n^2}(\phi_n(x)\phi_n'(y)-\phi_n(y)\phi_n'(x)+O(n^{-1/3+\varepsilon}))$$
Since the function is analytic, Lemma \ref{principe-maximum}
implies that this holds also for $x=y$.

As a consequence of Lemma \ref{approxt}, 
the asymptotic of Lemma \ref{softestimate} is not modified by
replacing $\tilde{s}_n$ by $s_n$ because $s_n-\tilde{s}_n=O(n^{-1})$.
Therefore, again by Lemmas \ref{approxt} and \ref{softestimate}, 
the above left hand side tends towards the Airy kernel.

\begin{theorem}\label{airy}
For any $\varepsilon >0$,
\begin{equation*}
\frac{1}{h_n}K_{n}^{a_n,b_n}(s_n+\frac{x}{h_n},s+\frac{y}{h_n})
=Ai(x,y)+O(n^{-1/3+\varepsilon})
\end{equation*}
where $s_n$ and  $h_n$ are defined above Equation \eqref{ashn}, and
the latter gives an estimate for $h_n$.
This asymptotic holds uniformly on any compact subset of $\mathbb{R}$.
\end{theorem}

From this, it is possible to establish a central limit type theorem
for the largest eigenvalues

\begin{proposition}\label{ultraspherical}
In the ultraspherical case $a_n=b_n$,
$$\left((\lambda_1-s_n)h_n,(\lambda_2-s_n)h_n,\ldots ,(\lambda_n-s_n)
h_n,0,0,\ldots \right)$$
converges in distribution towards the Airy ensemble as $n\rightarrow\infty$, in the sense
that any finite dimensional marginal converges in distribution.
\end{proposition}

\begin{proof}
We prove the convergence in distribution 
of $(\lambda_1-s_n)h_n$, the general statement
being obtained by the same standard methods.
Observe that for $v\geq u$,
\begin{equation*}
\begin{split}
P(\Lambda_n^{a_n,b_n}\cap [s_n+\frac{u}{h_n},s_n+\frac{v}{h_n}]\neq\emptyset)\leq
P(\Lambda_n^{a_n,b_n}\cap [s_n+\frac{u}{h_n},\infty )\neq\emptyset) \\
\leq P(\Lambda_n^{a_n,b_n}\cap [s_n+\frac{v}{h_n},\infty )\neq\emptyset)+
P(\Lambda_n^{a_n,b_n}\cap [s_n+\frac{u}{h_n},s_n+\frac{v}{h_n}]\neq\emptyset)
\end{split}
\end{equation*}
According to Proposition 6.4 of \cite{Le}, in the
ultraspherical case $a_n=b_n=\alpha n$, there exists a 
constant $C$ such that for every $0<\varepsilon\leq 1$ and $n\geq 1$,
\begin{equation*}
P(\lambda^n_{1}\geq s_n(1+\varepsilon ))\leq Ce^{-n\varepsilon^{3/2}/C}
\end{equation*}
This constant exists again in our more general framework (the constant $C$
of \cite{Le} can be chosen to depend explicitly and continuously on $\alpha$).

By Theorem \ref{airy} and dominated convegence, the right hand side of
Equation \eqref{deift} on the set
$[s_n+\frac{u}{h_n},s_n+\frac{v}{h_n}]$ converges as $n\rightarrow\infty$.
Additionally, when
$v\rightarrow\infty$, $P(\Lambda_n^{a_n,b_n}\cap [s_n+\frac{v}{h_n},\infty )\neq\emptyset)$
tends towards zero by the result of \cite{Le}.
\end{proof}

\begin{remark}
\begin{itemize}
\item
A direct analysis from the asymptotics or the modified Jacobi Kernel
for the above result escaped us because our results only hold only
uniformly on compact sets. A direct approach would be very interesting
but for the moment the result of \cite{Le} remains unavoidable.
\item
It would be interesting to check the result of \cite{Le} in the most general
($a_n\neq b_n$) case. This would result in the convergence after rescaling
towards an Airy ensemble in full generality.
\end{itemize}
\end{remark}

\subsection{Kernel at the hard edge.}

According to papers of Kuijlaars et al \cite{math-ph/0204006}
in which Riemann-Hilbert techniques are used, we expect 
the relevant kernels at the hard edge to involve Bessel functions.
We present a short and self contained algebraic solution
to this question.
We make the following
\begin{assumption}\label{hard}
The sequence $b_n$ is constant (we denote its value by the nonnegative
integer $b$), and  defining $\alpha_n=a_n/n$, we assume
$$0\leq \liminf_n \alpha_n <\infty$$
\end{assumption}
The {\it Bessel kernel\it} is defined by
\begin{equation}\label{defbessel}
F_b(u,v)=\frac{J_b(\sqrt{u})\sqrt{v}J_b'(\sqrt{v})-
J_b(\sqrt{v})\sqrt{u}J_b'(\sqrt{u})}{2(u-v)}
\end{equation}
where $J_b$ is the usual Bessel function of the first kind and of order $b$,
satisfying

\begin{equation}
J_b(z)=(z/2)^b\sum_{n=0}^{\infty}\frac{(-1)^n}{(n+b)!n!}(z/2)^{2n}
\end{equation}

\begin{theorem}\label{bessel}
For any $u,v\in\mathbb{R}_+$,
\begin{equation*}
\frac{1}{2n^2(1+\alpha_n)}
K_n^{a_n,b}(-1+\frac{u}{2n^2(1+\alpha_n)},-1+\frac{v}{2n^2(1+\alpha_n)})
=  F_b(u,v)+O(n^{-1})
\end{equation*}
and this estimate is uniform for any compact set in $\mathbb{R}_+$.
\end{theorem}

\begin{proof}
Observe that
\begin{equation*}
J_b'(z)=-J_{b+1}(z)+bJ_b(z)/z
\end{equation*}
therefore 
\begin{equation*}
F_b(u,v)=\frac{-J_b(\sqrt{u})\sqrt{v}J_{b+1}(\sqrt{v})+
J_b(\sqrt{v})\sqrt{u}J_{b+1}(\sqrt{u})}{2(u-v)}
\end{equation*}

For future use, note that the coefficient in
$u^{b/2+k}v^{b/2+l}$ of
$$-J_b(\sqrt{u})\sqrt{v}J_{b+1}(\sqrt{v})+ J_b(\sqrt{v})\sqrt{u}J_{b+1}(\sqrt{u})$$ 
is
\begin{equation}\label{coefft}
\frac{(k-l)}{(b+k)!k!(b+l)!l!}
\end{equation}

According to \cite{MR58:1288}, p. 7 Formula (2.2)
\begin{equation*}
P_n^{a,b}(-1+\frac{2u}{n^2})
=\frac{(b+1)_n}{n!}
\sum_{k=0}^n\frac{(-n)_k(n+a+b+1)_k}{k!(b+1)_k}
\left(\frac{u}{n^2}\right)^k
\end{equation*}
where $(a)_k=a(a+1)\ldots (a+k-1)$.

In order to establish asymptotic properties of this polynomial, let us write it in
the equivalent form
\begin{equation*}
P_n^{a,b}(-1+\frac{2u}{n^2})
=\frac{(b+n)!}{n!}
\sum_{k=0}^n\frac{(-n)_k(n+a+b+1)_k}{k!(b+k)!}
\left(\frac{u}{n^2}\right)^k
\end{equation*}
Isolate one generic summand of
$P_n^{a_n,b_n}(-1+\frac{u}{2n^2 (1+\alpha_n)})$:
\begin{equation*}
\frac{(b_n+n)!}{n!}
\frac{(-n)_k(n+a_n+b_n+1)_k}{k!(b_n+k)!}
\left(\frac{u}{4n^2 (1+\alpha_n)}\right)^k
\end{equation*}

The coefficient in $u^kv^l$ of 
\begin{eqnarray*}
P_n^{a_n,b_n}(-1+\frac{u}{n^2 2(1+\alpha_n)})P_{n-1}^{a_n,b_n}(-1+\frac{v}{2n^2 (1+\alpha_n)})\\
-P_n^{a_n,b_n}(-1+\frac{v}{n^2 2(1+\alpha_n)})P_{n-1}^{a_n,b_n}(-1+\frac{u}{2n^2 (1+\alpha_n)})
\end{eqnarray*}
is
\begin{eqnarray*}
\frac{(b+n)!(b+n-1)!}{n!(n-1)!k!(b+k)!l!(b+l)!(4n^2(1+\alpha_n))^{k+l}}\\
\{ (-n)_k(n+a_n+b+1)_k(-n+1)_l(n+a_n+b)_l-\\
(-n)_l(n+a_n+b+1)_l(-n+1)_k(n+a_n+b)_k\}\\
=\frac{(b+n)!(b+n-1)!}{n!(n-1)!k!(b+k)!l!(b+l)!(4n^2(1+\alpha_n))^{k+l}}\\
(-n+1)_{l-1}(n+a_n+b+1)_{l-1}(-n+1)_{k-1}(n+a_n+b+1)_{k-1}\\
\{-n(n+a_n+b+k)(-n+l)(n+a_n+b)-\\
(-n)(n+a_n+b+l)(-n+k)(n+a_n+b)\}
\end{eqnarray*}
This simplifies to
\begin{eqnarray*}
\frac{(b+n)!(b+n-1)!}{n!(n-1)!k!(b+k)!l!(b+l)!(4n^2(1+\alpha_n))^{k+l}}\\
(-n+1)_{l-1}(n+a_n+b+1)_{l-1}(-n+1)_{k-1}(n+a_n+b+1)_{k-1}\\
(-n)(n+a_n+b)(l-k)(2n+a_n+b)\\
=\frac{(b+n)!(b+n-1)!}{n!(n-1)!k!(b+k)!l!(b+l)!(4n^2(1+\alpha_n))^{k+l}}\\
(-n)_{l}(n+a_n+b)_{l}(-n+1)_{k-1}(n+a_n+b+1)_{k-1}
(l-k)(2n+a_n+b)
\end{eqnarray*}
As $n \to\infty$, the above expression can be simplified:
\begin{equation*}
\frac{(l-k)}{k!(b+k)!l!(b+l)!}n^{2b-2}(1+\alpha_n)^{-1}(-1)^{l+k-1}(2+\alpha_n)(1+0(n^{-1}))
\end{equation*}
One can show that the remainder function $0(n^{-1})$ is smaller than
$n^{-1}p(k,l)$ where $p(k,l)$ is a suitably chosen polynomial in $k$ and $l$).
Additionally,
 $\gamma_n=n^2(1+\alpha_n)^{b+1}(2+\alpha_n)^{-1}+O(n)$
and
\begin{equation*}
w^{a_n,b}(-1+\frac{u}{2n^2(1+\alpha_n)})=
2^{a_n}(\frac{u}{2n^2(1+\alpha_n)})^b(1+O(n^{-1}))
\end{equation*}
Together with Equation \eqref{coefft}, we deduce that
\begin{eqnarray*}
\frac{u-v}{2n^2(1+\alpha_n)}K_n^{a_n,b}(-1+\frac{u}{2n^2(1+\alpha_n)},-1+\frac{v}{2n^2(1+\alpha_n)})\\
=\frac{n!(n+a_n+b)!}{(n+a_n-1)!(n+b-1)!}\frac{2^{-a_n-b}}{2n+a_n+b}\\
w^{a_n/2,b/2}(-1+\frac{u}{2n^2(1+\alpha_n)})
w^{a_n/2,b/2}(-1+\frac{u}{2n^2(1+\alpha_n)})\\
(P_n^{a_n,b}(-1+\frac{u}{2n^2 (1+\alpha_n)})P_{n-1}^{a_n,b}(-1+\frac{v}{2n^2 (1+\alpha_n)})\\
-P_n^{a_n,b}(-1+\frac{v}{2n^2 (1+\alpha_n)})P_{n-1}^{a_n,b}(-1+\frac{u}{2n^2 (1+\alpha_n)}))\\
=(u-v)F_b(u,v)(1+O(n^{-1}))
\end{eqnarray*}

Therefore, for any compact subset of $\mathbb{C}\times\mathbb{C}$
not intersecting the 
diagonal $\{(x,x),x\in\mathbb{C}\}$, one has the announced 
result. By Lemma \ref{principe-maximum},
this convergence is 
uniform on any compact subset of $\mathbb{C}\times\mathbb{C}$.
\end{proof}

From this, we can also state a result about the behavior of the smallest
non-zero eigenvalue:

\begin{proposition}
Under Assumption \ref{hard}, the random sequence of vectors
$\left((\lambda_1+1)n^{2}(1+\alpha_n),(\lambda_2+1)n^{2}(1+\alpha_n),\ldots
,(\lambda_n+1)n^2(1+\alpha_n), \ldots \right)$
converges in distribution. 
\end{proposition}

\begin{proof}
This is a consequence of Equation \eqref{deift} and dominated convergence
theorem.
\end{proof}

\section{Concluding remarks.}

\subsection{Limiting procedures and modified Laguerre ensemble.}

In \cite{Co}, the following result is proved

\begin{theorem}\label{densite}
For $n\geq q$, $n\geq q'$,
let  $\pi_{n,q,q'}$ be the canonical projection of $\M{n}$ onto its 
upper left corner $\mathbb{M}_{q,q'}(\mathbb{C})$ with $q$ lines
and $q'$ columns. Let $dA$ be the standard Lebesgue measure
on $\mathbb{M}_{q,q'}(\mathbb{C})$
For $n\geq 2q \geq 2q'$, 
\begin{equation*}
\pi_{n,q,q'}^*(\mu_n)=c_{q,q',n}\det (1-AA^*)^{n-q-q'}1_{||A||\leq 1}dA
\end{equation*}
where $c_{q,q',n}$ is a normalization constant. 
\end{theorem}
This result is also a consequence of the present paper and it implies
\begin{theorem}\label{del}
Let $\nu_q$ be the probability measure
$c_{n,q}e^{-q\Tr MM^*}dM$ on $\M{q}$, and $q_n$ be a sequence of integers
tending towards infinity such that there exists a $C >0$
such that $q_n^{3}\leq C n$.
Then 
\begin{equation*}
|\sqrt{n/q_n}\pi_{n,q_n}^*(\mu_n)-\nu_{q_n}|=o(1)
\end{equation*}
where $|\cdot |$ denotes the total variation measure.
\end{theorem}

This result was already known to \cite{MR94g:60065} under the assumption that
$q_n^3=o(n)$.
Jiang informed the author that he recently obtained by different methods 
(\cite{jiang}) an improvement of this theorem to the case $q_n=o(n^2)$.

The Laguerre Polynomial $(L^{a}_n)_{n\geq 0}$ is a family of orthogonal
polynomials with respect to the measure $x^ae^{-x}1_{[0,\infty )}$ 
such that the leading coefficient is $(-1)^n/n!$.
These polynomials determine the determinantal point process
structure of a Wishart ensemble of parameter $(n,n+a,n^{-1})$.
As $n\rightarrow\infty$, $a_n/n\sim \alpha\geq 0$,
the average eigenvalues counting measure of $W(n,a_n,n^{-1})$
converges towards the so-called {\it Marcenko-Pastur\it} distribution 
$$const.\frac{\sqrt{(u-x)(x-v)}}{x}1_{[u,v]}dx$$
where $u=2+\alpha-2\sqrt{1+\alpha},v=2+\alpha+2\sqrt{1+\alpha}$

With Speicher's non-crossing cumulants theory (see \cite{NS}), one can prove that 
this distribution is both a free chi-square
distribution and a free Poisson distribution. 

Upon knowing that the average eigenvalue counting measure of the
$GUE$ converges towards the semi-circle distribution, it is easy to understand
via matrix models why the Marcenko-Pastur distribution is a free chi-square
distribution.

However, as far as the author knows, there was no matrix-model 
explanation for the coincidence between free Poisson distribution and 
Marcenko-Pastur distribution.
This paper provides an explanation: indeed, Theorem \ref{del}
this implies that contraction of a random projection by small projections 
is almost a Wishart matrix; in addition, contraction of a matrix $A$  
by a random projection of rank $d\alpha$
is a matrix model for $\alpha^{-1}$ fold free additive convolution of $A$.

Note that at the level of orthogonal polynomials, 
this is in accordance with the following 
well known approximation result in orthogonal polynomial theory
(see \cite{MR58:1288}, Formula (6.11)):
$$\lim_{b\rightarrow\infty}P_n^{a,b}(1-\frac{2x}{b})=L_n^a(x)$$

This gives a free probabilistic motivation for computing the local spacing
results for the Laguerre ensemble. The methods 
of this paper can me followed line by line to obtain similar results for
the modified Laguerre ensemble.  
For the asymptotic kernel at the edge we have to make use of the
differential equation $g''=\chi g$
with $g=e^{-x/2}x^{(a+1)/2}L_n^a(x)$ and
$$\chi (x)=\frac{2n+a+1}{2x}+\frac{1-a^2}{4x^2}-\frac{1}{4}$$ 
and at the hard edge, one can use the following hypergeometric representation
$$L_n^a(x)={n+a \choose a}M(-n,a+1,x)$$
where $M$ is the confluent hypergeometric function.

\subsection{Remaining questions.}

It would be very interesting to investigate possible extension of 
the dictionary between free probability and classical (possibly modified)
polynomials (Charlier, Meixner, etc...), and investigate universality properties.

The study of the semi group of free additive convolution, in particular
the study of random matrix models obtained by contractions of unitarily invariant
matrices, and obtaining nice Theorem \ref{jacobi}-like explicit densities
which can be handled for local asymptotics purposes remains a challenging 
problem in full generality, for which  other ideas are needed.

\bibliographystyle{alpha}
\bibliography{bibptrf}

\begin{thebibliography}{VDN92}

\bibitem[AS92]{MR94b:00012}
Milton Abramowitz and Irene~A. Stegun, editors.
\newblock {\em Handbook of mathematical functions with formulas, graphs, and
  mathematical tables}.
\newblock Dover Publications Inc., New York, 1992.
\newblock Reprint of the 1972 edition.

\bibitem[Ask75]{MR58:1288}
Richard Askey.
\newblock {\em Orthogonal polynomials and special functions}.
\newblock Society for Industrial and Applied Mathematics, Philadelphia, Pa.,
  1975.

\bibitem[BG99]{MR2001i:30040}
Christof Bosbach and Wolfgang Gawronski.
\newblock Strong asymptotics for {J}acobi polnomials with varying weights.
\newblock {\em Methods Appl. Anal.}, 6(1):39--54, 1999.
\newblock Dedicated to Richard A.\ Askey on the occasion of his 65th birthday,
  Part I.

\bibitem[CC02]{CC}
M.~Capitaine and M.~Casalis.
\newblock Asymptotic freeness by generalized moments for {G}aussian and
  {W}ishart matrices. {A}pplications to {B}eta random matrices.
\newblock {\em To appear in {I}ndiana {U}niversity {M}athematics {J}ournal},
  November 2002.

\bibitem[CI91]{MR92g:33011}
Li-Chen Chen and Mourad E.~H. Ismail.
\newblock On asymptotics of {J}acobi polynomials.
\newblock {\em SIAM J. Math. Anal.}, 22(5):1442--1449, 1991.

\bibitem[Col03a]{Co}
Beno\^\i{}t Collins.
\newblock Int\'egrales matricielles et probabilit\'es non-commutatives.
\newblock {\em Th\'ese de doctorat de l'Universit\'e Paris 6}, 2003.

\bibitem[Col03b]{imrn}
Benoit Collins.
\newblock Moments and cumulants of polynomial random variables on unitary
  groups, the {I}tzykson-{Z}uber integral and free probability.
\newblock {\em IMRN}, 17:953--982, 2003.

\bibitem[Dei99]{MR2000g:47048}
P.~A. Deift.
\newblock {\em Orthogonal polynomials and random matrices: a
  {R}iemann-{H}ilbert approach}, volume~3 of {\em Courant Lecture Notes in
  Mathematics}.
\newblock New York University Courant Institute of Mathematical Sciences, New
  York, 1999.

\bibitem[DEL92]{MR94g:60065}
Persi~W. Diaconis, Morris~L. Eaton, and Steffen~L. Lauritzen.
\newblock Finite de {F}inetti theorems in linear models and multivariate
  analysis.
\newblock {\em Scand. J. Statist.}, 19(4):289--315, 1992.

\bibitem[Dou03]{Do}
Yan Doumerc.
\newblock Matrix {J}acobi process.
\newblock {\em Work in progress}, 2003.

\bibitem[For02]{Fo}
Peter Forrester.
\newblock {\em Log-gases and Random matrices, Chapter 2}.
\newblock http://www.ms.unimelb.edu.au/~matpjf/matpjf.html, 2002.

\bibitem[GS91]{MR92j:33019}
Wolfgang Gawronski and Bruce Shawyer.
\newblock Strong asymptotics and the limit distribution of the zeros of
  {J}acobi polynomials {$P\sb n\sp {(an+\alpha,bn+\beta)}$}.
\newblock In {\em Progress in approximation theory}, pages 379--404. Academic
  Press, Boston, MA, 1991.

\bibitem[Jia03]{jiang}
Tiefeng Jiang.
\newblock Maxima of entries of {H}aar distributed matrices.
\newblock {\em preprint, available at
  http://www.stat.umn.edu/~tjiang/papers/haar1.pdf}, 2003.

\bibitem[Joh01]{MR2002j:15024}
Kurt Johansson.
\newblock Universality of the local spacing distribution in certain ensembles
  of {H}ermitian {W}igner matrices.
\newblock {\em Comm. Math. Phys.}, 215(3):683--705, 2001.

\bibitem[KV02]{math-ph/0204006}
Arno Kuijlaars and Maarten Vanlessen.
\newblock Universality for eigenvalue correlations from the modified {J}acobi
  unitary ensemble.
\newblock {\em IMRN}, 30:1575--1600, 2002.

\bibitem[Led02]{Le}
Michel Ledoux.
\newblock Differential operators and spectral distributions of invariant
  ensembles from the classical orthogonal polynomials part {I}: the continuous
  case.
\newblock {\em To appear in {E}lect. {J}. {P}robab.}, November 2002.

\bibitem[Meh91]{MR92f:82002}
Madan~Lal Mehta.
\newblock {\em Random matrices}.
\newblock Academic Press Inc., Boston, MA, second edition, 1991.

\bibitem[MSV79]{MR80g:33021}
D.~S. Moak, E.~B. Saff, and R.~S. Varga.
\newblock On the zeros of {J}acobi polynomials {$P\sb{n}\sp{(\alpha
  \sb{n},\beta \sb{n})}(x)$}.
\newblock {\em Trans. Amer. Math. Soc.}, 249(1):159--162, 1979.

\bibitem[NS00]{NS}
A.~Nica and R.~Speicher.
\newblock Lectures notes of the free probability semester at {I}{H}{P}.
\newblock 2000.

\bibitem[Ol'90]{0724.22020}
G.I. Ol'shanskij.
\newblock Unitary representations of infinite dimensional pairs (g,k) and the
  formalism of {R}. {H}owe.
\newblock {\em Representation of Lie groups and related topics, Adv. Stud.
  Contemp. Math. 7}, pages 269--463, 1990.

\bibitem[Sos99]{MR2001i:82037}
Alexander Soshnikov.
\newblock Universality at the edge of the spectrum in {W}igner random matrices.
\newblock {\em Comm. Math. Phys.}, 207(3):697--733, 1999.

\bibitem[Sze75]{MR51:8724}
G{\'a}bor Szeg{\H{o}}.
\newblock {\em Orthogonal polynomials}.
\newblock American Mathematical Society, Providence, R.I., fourth edition,
  1975.
\newblock American Mathematical Society, Colloquium Publications, Vol. XXIII.

\bibitem[VDN92]{MR94c:46133}
D.~V. Voiculescu, K.~J. Dykema, and A.~Nica.
\newblock {\em Free random variables}.
\newblock American Mathematical Society, Providence, RI, 1992.
\newblock A noncommutative probability approach to free products with
  applications to random matrices, operator algebras and harmonic analysis on
  free groups.

\bibitem[Voi98]{MR2000d:46080}
Dan Voiculescu.
\newblock A strengthened asymptotic freeness result for random matrices with
  applications to free entropy.
\newblock {\em Internat. Math. Res. Notices}, (1):41--63, 1998.

\bibitem[Xu97]{MR99f:81185}
Feng Xu.
\newblock A random matrix model from two-dimensional {Y}ang-{M}ills theory.
\newblock {\em Comm. Math. Phys.}, 190(2):287--307, 1997.

\end{thebibliography}

\end{document}